# IMPROVED APPROXIMATION ALGORITHM FOR SOLUTION OF NONLINEAR PARTIAL DIFFERENTIAL EQUATIONS

Sema SERVI[a] , Yildiray KESKIN[b] , Galip OTURANC[b]

[a]*Vocational School of Tecnical Sciences*,

[b]*Department of Mathematics*, *Science Faculty*,

Selcuk University, Konya 42075, Turkey

**Abstract**

In this study,a new method was presented by developing Reduced differential transform method in order to find approximate solution of partial differential equations.

Here, RDTM with fixed grid size algorithm was developed for the first time for Reduced Differential Transform Method by dividing solution intervals given to us into fixed grids. The efficiency and advantage of this method was given in homogenous heat equation existing in literature and in the application part on Burger's equation. When approximate solution obtained by this new method and known exact solutions were compared, it is seen that there is definite consistence between both two solutions.

*AMS 2015: 65M55*

## 1. Introduction

The partial differential equation has been used in many mathematical, engineering problems,fluid mechanic, thermodynamic, heat transfer, physics and mathematical physics. In recent years with the development of the tecnology, many effective methods have been proposed for solving the partial differential equations such as variational iteration method [7], the homotopy analysis method [8],homotopy perturbation method[10], the differential transform method (DTM[11],the Adomian's decomposition method (ADM) [12-7],sine–cosine method,[13],reduced differential transform method (RDTM)[5,6,15,16,17,18] and so forth.

In this study, we developed RDTM with fixed grid size by making improvements on RDTM by these methods.

RDTM method is a method found out by Y. Keskin in 2009 with the development of DTM and it was applied to a variety of problems.[1,6,15-18]

DTM was introduced by Zhou in 1986 for the first time for the solution of linear and non-linear initial value problem that was faced in the electric and electric circuit analysis. With this method, complicated integral statements faced in the equations were eliminated by transforming partial differentiable equations into algebraic equations. Also results were obtained by simple operations. However, RDTM was developed to get results more quickly by reducing operations faced while calculating approximate result and reducing the number of iterations. When this method was compared to other known methods, it is seen that RDTM is more efficient and more effective than others. [5,6,15,17].

It is obvious that, there is a need for improvement in present numerical methods to get the approximate solution faster or more sensitive and efficient results. From here, we developed RDTM with fixed grid size algorithm by dividing the intervals given in the equations into equal parts. Fixed grid size algorithms was earlier applied by Ming-Jyi Jang , Chieh-Li Chen , Yung-Chin Liy for the approximate solution of linear and non-linear initial value problems with DTM and positive results were obtained. We compared the results obtained with known exact solution by applying RDTM which is more efficient than DTM . To better understand , we also applied RDTM on the heat and Burger equations.

2. **RDTM with grid size solution**

The basic definitions of reduced differential transform method are introduced as in the references[1,2,3,4,5]:

**Definition 2.1.** If the function $u(x,t)$ is analytic and differentiated continuously with respect to time t and space x in the domain of interest, the equation where the t-dimensional spectrum function $U_k(x)$ is the transformed function is like this:

$$U_k(x) = \frac{1}{k!}\left[\frac{\partial^k}{\partial t^k}u(x,t)\right]_{t=0} \quad (2.1)$$

. In this paper, the lowercase $u(x,t)$ represents the original function while the uppercase $U_k(x)$ stands for the transformed function.

The differential inverse transform of $U_k(x)$ is defined as in the following:

$$u(x,t) = \sum_{k=0}^{\infty} U_k(x)t^k. \quad (2.2)$$

When equations of (2.1) and (2.2) are combined, we write

$$u(x,t) = \sum_{k=0}^{\infty} \frac{1}{k!}\left[\frac{\partial^k}{\partial t^k}u(x,t)\right]_{t=0}. \quad (2.3)$$

From the above definitions, it can be found that the concept of the reduced differential transform is derived from the power series expansion.

According to the RDTM and Table 1, we can construct the following iteration formula:

**Table 1. Reduced differential transformation**

| Functional Form | Transformed Form |
|---|---|
| $u(x,t)$ | $U_k(x) = \dfrac{1}{k!}\left[\dfrac{\partial^k}{\partial t^k}u(x,t)\right]_{t=0}$ |
| $w(x,t) = u(x,t) \pm v(x,t)$ | $W_k(x) = U_k(x) \pm V_k(x)$ |
| $w(x,t) = \alpha u(x,t)$ | $W_k(x) = \alpha U_k(x)$ ($\alpha$ is a constant) |
| $w(x,y) = x^m t^n$ | $W_k(x) = x^m \delta(k-n)$ |
| $w(x,y) = x^m t^n u(x,t)$ | $W_k(x) = x^m U(k-n)$ |
| $w(x,t) = u(x,t)v(x,t)$ | $W_k(x) = \sum_{r=0}^{k} V_r(x) U_{k-r}(x) = \sum_{r=0}^{k} U_r(x) V_{k-r}(x)$ |
| $w(x,t) = \dfrac{\partial^r}{\partial t^r} u(x,t)$ | $W_k(x) = (k+1)...(k+r) U_{k+1}(x) = \dfrac{(k+r)!}{k!} U_{k+r}(x)$ |
| $w(x,t) = \dfrac{\partial}{\partial x} u(x,t)$ | $W_k(x) = \dfrac{\partial}{\partial x} U_k(x)$ |

Our purpose here is to investigate the solution of (1.2) the interval $[0;T]$ and the interval is divided into $N$ equally spaced subintervals with th egridpoints $\{t_0, t_1, t_2, ..., t_N\}$ as indicated in Figure 1

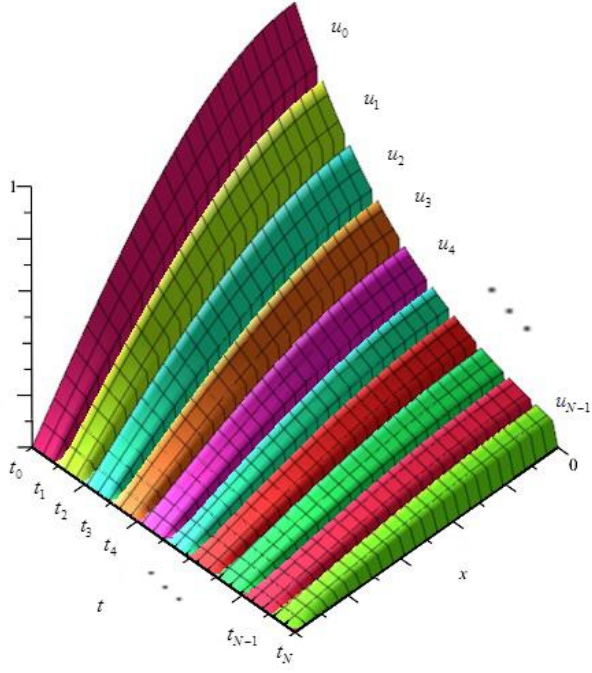

**Figure1: Approximation functions in each sub domain**

Here $t_0 = 0, t_N = T$ as in the formulation

$$t_i = t_0 + ih, \text{ for each i=0,1,2,...,N and } h = \frac{T}{N} \qquad (2.4)$$

Here approximate solution function $u(x,t)$ can be indicated as $u^0(x,t)$ at the first interval and after RDTM procedures are applied, Taylor polynomials from nth order can be obtained at the point $t_0 = 0$ as

$$u^0(x,t) = U_0^0(x,t_0) + U_1^0(x,t_0)(t-t_0) + U_2^0(x,t_0)(t-t_0)^2 + ... + U_n^0(x,t_0)(t-t_0)^n \qquad (2.5)$$

Fromtheinitialcondition (1.2), we can write,

$$u(x,t_0) = U_0^0(x,t_0) \qquad (2.6)$$

Therefore the approximate value of the function $u(x,t)$ at the point $t_1$ can be evaluatedby (2.5) as

$$\begin{aligned}
u(x,t_1) &\approx u^0(x,t_1) \\
&= U_0^0(x,t_0) + U_1^0(x,t_0)(t_1-t_0) + U_2^0(x,t_0)(t_1-t_0)^2 + ... + U_n^0(x,t_0)(t_1-t_0)^n \\
&= U_0^0 + U_1^0 h + U_2^0 h^2 + ... + U_n^0 h^n \\
&= \sum_{j=0}^{n} U_j^0 h^j
\end{aligned} \qquad (2.7)$$

Since the initial value of the second sub-domain is equal to the value of $u^0(x,t)$ at the point $t_1$ we can write

$$u^1(x,t) = U_0^1(x,t_1) + U_1^1(x,t_1)(t-t_1) + U_2^1(x,t_1)(t-t_1)^2 + ... + U_n^1(x,t_1)(t-t_1)^n$$

$$u^1(x,t_1) = u^0(x,t_1) = U_0^1(x,t_1) \qquad (2.8)$$

In a similar way, $u^2(x,t_2)$ can be calculated as

$$u(x,t_2) \approx u^1(x,t_2)$$
$$= U_0^1(x,t_1) + U_1^1(x,t_1)(t_2 - t_1) + U_2^1(x,t_1)(t_2 - t_1)^2 + \ldots + U_n^1(x,t_1)(t_2 - t_1)^n$$
$$= U_0^1 + U_1^1 h + U_2^1 h^2 + \ldots + U_n^1 h^n$$
$$= \sum_{j=0}^{n} U_j^1 h^j \tag{2.9}$$

Once following the same procedure, we can obtain the solution of $u^i(x,t)$ at the grid point $t_{i+1}$ as

$$u(x,t_{i+1}) \approx u^i(x,t_{i+1})$$
$$= U_0^i(x,t_i) + U_1^i(x,t_i)(t_{i+1} - t_i) + U_2^i(x,t_i)(t_{i+1} - t_i)^2 + \ldots + U_n^i(x,t_i)(t_{i+1} - t_i)^n$$
$$= U_0^i + U_1^i h + U_2^i h^2 + \ldots + U_n^i h^n$$
$$= \sum_{j=0}^{n} U_j^i h^j \tag{2.10}$$

From here, if the values of $u^i(x,t)$ are calculated, analytical solution of $u(x,t)$ is obtained as

$$u(x,t) = \begin{cases} u^0(x,t), & 0 \leq x \leq 1, \ 0 < t \leq t_1 \\ u^1(x,t), & 0 \leq x \leq 1, \ t_1 \leq t \leq t_2 \\ \vdots \\ u^{N-1}(x,t), & 0 \leq x \leq 1, \ t_{N-1} \leq t \leq t_N \end{cases}$$

We get more sensitive results when we increase the value of $N$.

### 3. Application

In this section, in order to assess efficiency and advantages of our method the accuracy of the method was compared with the exact solutions.[1]. The obtained numerical results are very encouraging. We consider the following examples.

**Example1:** Homogeneous heat equation is

$$u_t = u_{xx}, \quad 0 < x < 1, 0 < t \leq T \tag{3.1}$$

the initial conditions are given

$$u(x,0) = f(x) = \cos(\frac{\pi}{2}x), \ 0 < x < 1 \tag{3.2}$$

And the exact solution for the (3.1) equation is

$$u(x,t) = \exp(-\frac{\pi^2}{4}t)\cos(\frac{\pi}{2}x)$$

By choosing N=5 we get t values as $\left\{ t_0 = 0, t_1 = \frac{1}{5}, t_2 = \frac{2}{5}, t_3 = \frac{3}{5}, t_4 = \frac{4}{5}, t_5 = 1 \right\}$

For the solution procedure, we first take the reduced differential transform of (3.1) and (3.2) by the aid of Table 1 and we get the following equation

$$(k+1)U^0_{k+1}(x) = \frac{\partial^2}{\partial x^2}U^0_k(x), \quad u(x,t_0) = U^0_0(x,t_0) = \cos(\frac{\pi}{2}x)$$

The transformation $U^0_k(x)$ are indicated at a few steps for k values and $t_0 = 0$

are as in the following

for k=0,

$$U^0_1(x) = \frac{\partial^2}{\partial x^2}U^0_0(x)$$

$$U^0_1(x) = \frac{\partial^2}{\partial x^2}\cos(\frac{\pi}{2}x) = -\frac{1}{4}\cos(\frac{\pi}{2}x)\pi^2 = 2.467401101\cos(1.570796327.x)$$

for k=1,

$$U^0_2(x) = \frac{1}{2}\frac{\partial^2}{\partial x^2}U^0_1(x) = \frac{1}{32}\cos(\frac{\pi}{2}x)\pi^4 = 3.044034097\cos(1.570796327.x)$$

$$\vdots$$

If the values of $N=5$ and $U^0_n(x,t_0)$ are written in the places in the following equations, the approximate solution can be obtained as $u^0(x,t)$

$$u^0(x,t) = \cos(\frac{\pi}{2}x) - 2.467401101\cos(1.570796327x)t + 3.044034097\cos(1.570796327x)t^2$$
$$- 2.503617694\cos(1.570796327x)t^3 + 1.544357264\cos(1.570796327x)t^4$$
$$- 0.7621097628\cos(1.570796327x)t^5$$

for $t_1 = \frac{1}{5}$ the equation and initial condition for (3.1) is

$$(k+1)U^1_{k+1}(x) = \frac{\partial^2}{\partial x^2}U^1_k(x), \quad u(x,t_1) = u^0(x,t_1) = 0.6104792987\cos(1.570796327x)$$

For k values and $t_0 = 0$ here, the transformation $U^0_k(x)$ as indicated:

$$u^1(x,t) = 0.6104792987\cos(1.570796327x) - 1.506297294\cos(1.570796327x)(-1/5+t)$$
$$+1.858319801\cos(1.570796327x)(-1/5+t)^2 - 1.528406774\cos(1.570796327x)(-1/5+t)^3$$
$$+0.9427981392\cos(1.570796327x)(-1/5+t)^4 - 0.4652522332\cos(1.570796327x)(-1/5+t)^5$$

For $t_2 = \frac{2}{5}$, the equation of (3.1) and the initial condition are,

$$(k+1)U^2_{k+1}(x) = \frac{\partial^2}{\partial x^2}U^2_k(x), \quad u(x,t_2) = u^1(x,t_2) = 0.3726849740\cos(1.570796327x)$$

$$u^2(x,t) = 0.3726849740\cos(1.570796327x) - 0.9195633152\cos(1.570796327x)(-2/5+t)$$
$$+ 1.134465768\cos(1.570796327x)(-2/5+t)^2 - 0.9330606947\cos(1.570796327x)(-2/5+t)^3$$
$$+ 0.5755587462\cos(1.570796327x)(-2/5+t)^4 - 0.2840268568\cos(1.570796327x)(-2/5+t)^5$$

For $t_3 = \dfrac{3}{5}$ the equation of (3.1) and the initial condition are,

$$(k+1)U^3_{k+1}(x) = \frac{\partial^2}{\partial x^2}U^3_k(x), \quad u(x,t_3) = u^2(x,t_3) = 0.2275164615\cos(1.570796327x)$$

$$u^3(x,t) = 0.2275164615\cos(1.570796327x) - 0.5613743677\cos(1.570796327x)(-3/5+t)$$
$$+ 0.6925678665\cos(1.570796327x)(-3/5+t)^2 - 0.5696142387\cos(1.570796327x)(-3/5+t)^3$$
$$+ 0.3513667000\cos(1.570796327x)(-3/5+t)^4 - 0.1733925165\cos(1.570796327x)(-3/5+t)^5$$

For $t_4 = \dfrac{4}{5}$, the equation of (3.1) and the initial condition are,

$$(k+1)U^4_{k+1}(x) = \frac{\partial^2}{\partial x^2}U^4_k(x), \quad u(x,t_4) = u^3(x,t_4) = 0.1388940899\cos(1.570796327x)$$

$$u^4(x,t) = 0.1388940899\cos(1.570796327x) - 0.3427074304\cos(1.570796327x)(-4/5+t)$$
$$+ 0.4227983455\cos(1.570796327x)(-4/5+t)^2 - 0.3477377010\cos(1.570796327x)(-4/5+t)^3$$
$$+ 0.2145020966\cos(1.570796327x)(-4/5+t)^4 - 0.1058525419\cos(1.570796327x)(-4/5+t)^5$$

Furthermore, the absolute errors according to values of $x$ and $t$ are given as in the following tables.

**Table 2. Absolute Errors for solution of Example1**

|  | t | AbsoluteErrorfor N=4 | AbsoluteErrorfor N=16 |
|---|---|---|---|
| x=0,25 | 0,25 | 0.00006491050816721604 | $4.261649214506.10^{-8}$ |
|  | 0,5 | 0.00004968038314752197 | $3.264827855296.10^{-8}$ |
|  | 0,75 | 0.00002688688512066903 | $1.772463571363.10^{-8}$ |
|  | 1 | $1.1066788288112711837.10^{-10}$ | $1.1068228370224933920.10^{-10}$ |
| x=0,5 | 0,25 | 0.00007005226033954058 | $4.584309270421.10^{-8}$ |
|  | 0,5 | 0.00005361569560274446 | $3.510352789731.10^{-8}$ |
|  | 0,75 | 0.00002901663531705991 | $1.902785042357.10^{-8}$ |
|  | 1 | $5.9713205373261158297.10^{-11}$ | $5.9728746958876876434.10^{-11}$ |
| x=0,75 | 0,25 | 0.00005670114017469839 | $3.720636020014.10^{-8}$ |
|  | 0,5 | 0.00004339718291423823 | $2.848555088372.10^{-8}$ |
|  | 0,75 | 0.00002348640231079870 | $1.5432435640069.10^{-8}$ |
|  | 1 | $3.2219527463403285982.10^{-11}$ | $3.2232106987014619021.10^{-11}$ |
| x=1 | 0,25 | 0.00004079503870330856 | $2.6717947869064.10^{-8}$ |
|  | 0,5 | 0.00003122317574236045 | $2.0453907992587.10^{-8}$ |
|  | 0,75 | 0.00001689786966556789 | $1.1078299400940.10^{-8}$ |
|  | 1 | $1.7384729970441294615.10^{-11}$ | $1.7393780631399927727.10^{-11}$ |

**Then, following graphics show the behavior of the results in the above table.**

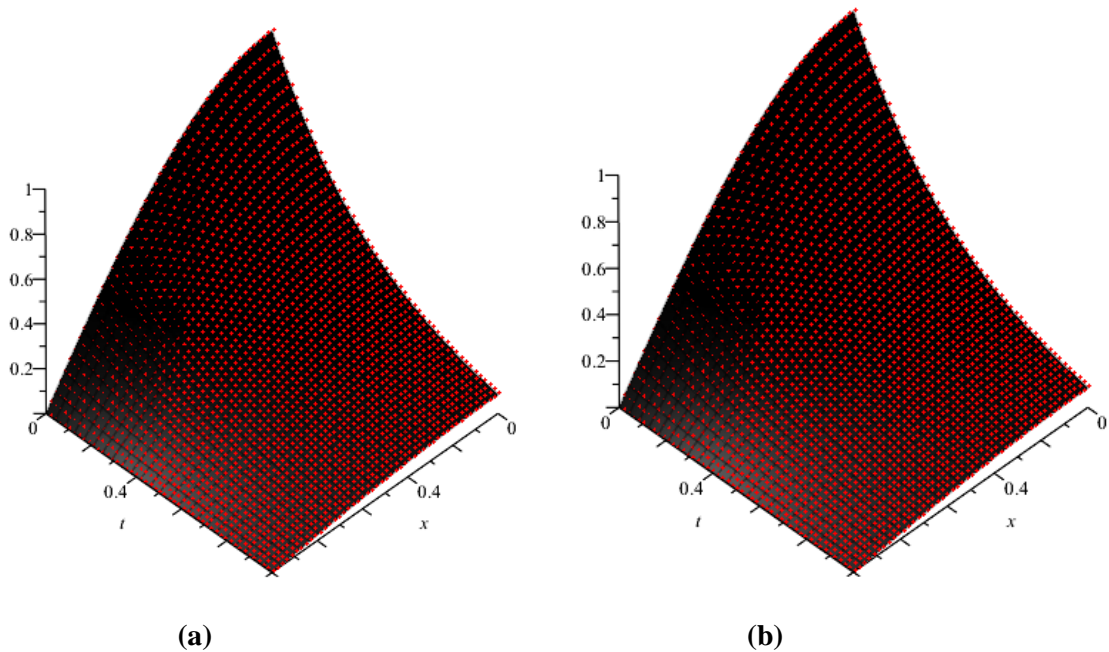

(a)　　　　　　　　　　　　　　　　　(b)

**Figure 2:In this graphics, the black area indicates the exact solution and the red spots indicate the result we obtained by our method. (a): for N=4, (b): for N=16. According to the figures, it is seen that there are a pretty good consistence between two solutions.**

**Example 2**

The Burgers' equation [14] is

$$u_t + uu_x - u_{xx} = 0, \quad 0 < x < 1, 0 < t \leq T \tag{3.3}$$

the initial conditions are given in (3.3)

$$u(x,0) = x \tag{3.4}$$

And the exact solution for the (3.3) equation is

$$u(x,0) = \frac{x}{1+t}, \quad |t| < 1$$

By choosing N=5 we get t values as $\left\{ t_0 = 0, t_1 = \frac{1}{5}, t_2 = \frac{2}{5}, t_3 = \frac{3}{5}, t_4 = \frac{4}{5}, t_5 = 1 \right\}$

For the solution procedure, we first take the reduced differential transform of (3.3) and (3.4) by the aid of Table 1 and we can get the following equation

$$(k+1)U^0_{k+1}(x) = -\sum_{r=0}^{k} U_r(x)\frac{\partial}{\partial x}U^0_{k-r}(x) + \frac{\partial^2}{\partial x^2}U^0_k(x), \quad u(x,t_0) = U^0_0(x,t_0) = x \quad (3.5)$$

Substituting (3.4) into (3.5), we obtain the following $U_k^{N-1}(x)$ values successively. Then we can write respectively approximate solution $u^{N-1}(x,t)$ for N=5 by calculating value of $U_k^{N-1}(x)$

$$u^0(x,t) = x - xt + xt^2 - xt^3 + xt^4 - xt^5$$

$$u^1(x,t) = 0.9999017774x - 0.9991822780xt + 0.9910716932xt^2 - 0.9374365902xt^3$$
$$+ 0.7365183841xt^4 - 0.3347693965xt^5$$

$$u^2(x,t) = 0.9993495237x - 0.9910875406xt + 0.9400113231xt^2 - 0.7700096803xt^3$$
$$+ 0.4513625924xt^4 - 0.1327495336xt^5$$

$$u^3(x,t) = 0.9971002908x - 0.9691870106xt + 0.8532651503xt^2 - 0.5958175754xt^3$$
$$+ 0.2740647794xt^4 - 0.05957774050xt^5$$

$$u^4(x,t) = 0.9921696552x - .934270172xt + 0.7536099398xt^2 - 0.452611796xt^3$$
$$+ 0.1704583443xt^4 - 0.2938872486xt^5$$

The following graphic shows the results of exact and approximate solution.

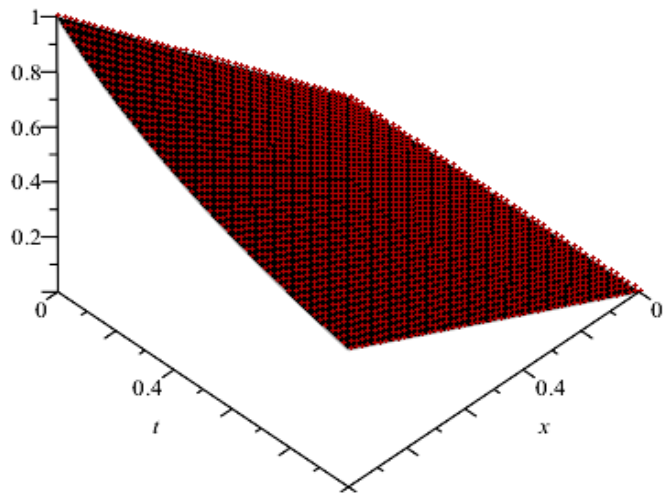

**Figure3: In this graphic, the black area indicates the exact solution for N=4 and the red spots indicate the result we obtained with our method.**

Furthermore, the exact and approximate solution according to values of x and t are given in the following Table 3 and Table 4

**Table3:Numerical Results of RDTM solution of Example2 for N=4**

|        | T    | Approximate solution      | Exact Solution            |
|--------|------|---------------------------|---------------------------|
| x=0,25 | 0,25 | 0.20000004750000000000    | 0.20000000000000000000    |
|        | 0,5  | 0.40000009500000000000    | 0.40000000000000000000    |
|        | 0,75 | 0.60000014250000000000    | 0.60000000000000000000    |
|        | 1    | 0.80000019000000000000    | 0.80000000000000000000    |
| x=0,5  | 0,25 | 0.16666670305000000000    | 0.16666666666666666667    |
|        | 0,5  | 0.33333340610000000000    | 0.33333333333333333333    |
|        | 0,75 | 0.50000010915000000000    | 0.50000000000000000000    |
|        | 1    | 0.80000019000000000000    | 0.80000000000000000000    |
| x=0,75 | 0,25 | 0.20000004750000000000    | 0.20000000000000000000    |
|        | 0,5  | 0.40000009500000000000    | 0.40000000000000000000    |
|        | 0,75 | 0.60000014250000000000    | 0.60000000000000000000    |
|        | 1    | 0.66666681220000000000    | 0.66666666666666666667    |
| x=1    | 0,25 | 0.12500002082541240453    | 0.12500002082541240453    |
|        | 0,5  | 0.25000004165082480908    | 0.25000000000000000000    |
|        | 0,75 | 0.37500006247623721361    | 0.37500000000000000000    |
|        | 1    | 0.50000008330164961816    | 0.50000000000000000000    |

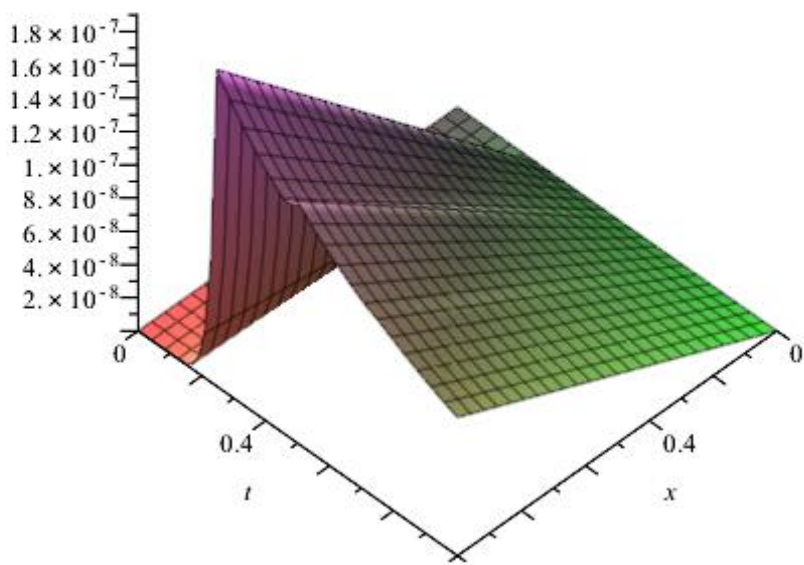

**Figure4: The surface shows the error of $u(x,t)$ for $N = 4$**

**Table4:Numerical Results of RDTM solution of Example2 for N=16**

|  | t | Approximatesolution | Exact Solution |
|---|---|---|---|
| x=0,25 | 0,25 | 0.19999999985000000000 | 0.20000000000000000000 |
|  | 0,5 | 0.39999999970000000000 | 0.40000000000000000000 |
|  | 0,75 | 0.59999999955000000000 | 0.60000000000000000000 |
|  | 1 | 0.79999999940000000000 | 0.80000000000000000000 |
| x=0,5 | 0,25 | 0.16666666657500000000 | 0.16666666666666666667 |
|  | 0,5 | 0.33333333315000000000 | 0.33333333333333333333 |
|  | 0,75 | 0.49999999972500000000 | 0.50000000000000000000 |
|  | 1 | 0.66666666630000000000 | 0.66666666666666666667 |
| x=0,75 | 0,25 | 0.14285714280000000000 | 0.14285714285714285714 |
|  | 0,5 | 0.28571428560000000000 | 0.28571428571428571429 |
|  | 0,75 | 0.42857142840000000000 | 0.42857142857142857143 |
|  | 1 | 0.57142857120000000000 | 0.57142857142857142857 |
| x=1 | 0,25 | 0.12499999998609105881 | 0.12500000000000000000 |
|  | 0,5 | 0.24999999997218211762 | 0.25000000000000000000 |
|  | 0,75 | 0.37499999995827317645 | 0.37500000000000000000 |
|  | 1 | 0.49999999994436423526 | 0.50000000000000000000 |

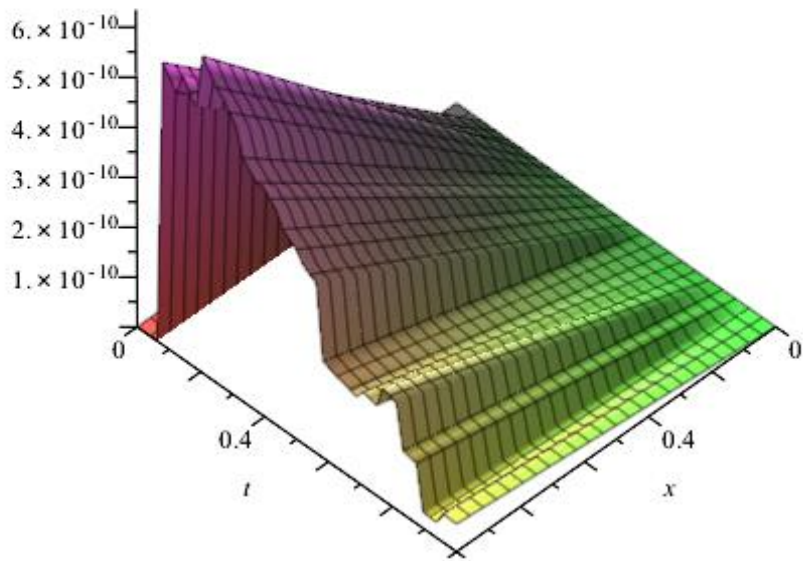

**Figure5: The surface shows the error** $u(x,t)$ **for** $N=16$

## 4. Conclusion

The reduced differential transform method with fixed grid size was first applied to the homogeneous heat equation and Burgers equation which are partial differential equation. The fixed grid algorithm is simply adaptable, sufficient and easily-programmable. Also, the used adaptive fixed grid size technique provides short and effective correction for the approximate solution and reduces the error. Two examples are given in order to show that the reduced differentials transform method with fixed grid size is a powerful mathematical tool for solving this partial differential equation. The main advantage of the method is the fact that it provides its users an analytical approximation in many cases. The exact solution in a rapidly convergent sequence is computed elegantly. In this study, for calculating the series obtained from the reduced differential transform method, we used the Maple programme.

## 5.References